\documentclass{article}
\usepackage[utf8]{inputenc}

\usepackage{graphicx, amsthm, amsmath, amssymb}
\newtheorem{theorem}{Theorem}

\newtheorem{lemma}[theorem]{Lemma}

\newtheorem{proposition}{Proposition}
\newtheorem{example}{Example}

\title{Rainbow triangles}
\author{Steven Senger}
\date{February 9, 2017}

\begin{document}

\maketitle

\begin{abstract}
We seek conditions under which colorings of various vector spaces are guaranteed to have a copy of a unit equilateral triangle, having each vertex in a different color class. In particular, we explore the analogous question in the setting of vector spaces finite fields, with an appropriate notion of distance.
\end{abstract}
\section{Introduction}

Classical Ramsey problems typically involve breaking a given ambient set up into disjoint subsets called color classes, then looking for conditions under which a given configuration will be present in one of the color classes. Here, we consider so-called anti-Ramsey or {\it rainbow} problems. Rainbow problems have been studied in different contexts; see \cite{AF} and \cite{JLMNR} and the references therein for some arithmetic rainbow results, and see \cite{Bab85} and \cite{LR} and the references therein for results in graph theory. Here we consider unit equilateral triangles whose vertices come from distinct color classes. The author would like to thank Janos Pach and Hans Parshall for helpful suggestions greatly improving the content of this paper.
\subsection{Background Information}
One of the largest motivating questions for this is the unit distance problem, due to Erd\H os, in \cite{Erd46}. In short, given a large finite number, $n$, what is the maximum number of times that two points can exactly a unit distance apart in a set of $n$ points in the plane? The current record of $cn^\frac{4}{3}$ was achieved by Spencer, Szemer\' edi, and Trotter in \cite{SST84}. People have studied similar questions involving triangles and higher-order simplices. See the wonderful survey book by Brass, Moser, and Pach, \cite{BMP00}.

The plane is not the only setting where such questions can be considered however. In \cite{IR07}, Iosevich and Rudnev considered distance problems in vector spaces over finite fields. Of course, they used an analog of distance, as vector spaces over finite fields do not admit Euclidean distances as such. There has been much activity in this area as people have considered related problems involving simplices, in \cite{HIKSU}, specifically triangles, in \cite{BIP}, and other functionals, such as dot products, in \cite{CHIKR}.

Specifically their initial result states that if one has a large enough subset of a vector space over a finite field, then there must be two points in that subset that are a given distance apart. We can think of this as a two-point configuration where the two points come from the same subset. In what follows, we will be looking for three-point configurations, where the three points all come from distinct sets.

Suppose that $q$ is a large, odd, prime power. Consider the two-dimensional vector space over the finite field with $q$ elements, $\mathbb F_q^2.$ In this setting, a common analog of Euclidean distance is the following functional, $d:\mathbb F_q^2\times\mathbb F_q^2\rightarrow \mathbb F_q,$ defined by 
$$d(x,y)=(x_1-y_1)^2+(x_2-y_2)^2,$$
which we shall, by abuse of language, call a ``distance." See \cite{IR07}, by Iosevich and Rudnev, for good introduction to this notion of distance.

 Also, given two functions, $f$ and $g$, we use the symbol $f \lesssim g$ to denote that $f=O(g).$ Similarly, we use $f\approx g$ if $f \lesssim g$ and $f \gtrsim g.$ Occasionally, we will use the letter $c$ to represent a constant. This should be clear by context. Also, we will frequently use the notation $f\not\lesssim 1$ to mean that $f$ grows faster than a constant with respect to the parameter.

\subsection{Main results}
Note that regardless of whether we are in $\mathbb R^2$ or $\mathbb F_q^2$, it is easy to see that if each point has a different color, every equilateral triangle will be rainbow. We first prove a simple result, Proposition \ref{size2}, with straightforward combinatorial methods as an introduction to the main result, Theorem \ref{main}.

\begin{proposition}\label{size2}
Given a coloring of $\mathbb R^2$ or $\mathbb F_q^2$ where no color class has more than two points, if there exist unit equilateral triangles, there must be a rainbow unit equilateral triangle.
\end{proposition}

The main result should be seen as an asymptotic result in $q$. That is, suppose we have a family of colorings of $\mathbb F_{q_j}$, for some sequence of $q_j$ tending to infinity. Then the constant, $c$, mentioned below should be independent of $q_j$.

\begin{theorem}\label{main}
Given a coloring of $\mathbb F_q^2$, where no color class has size greater than $cq^2,$ for any positive constant, $c$, if there exist unit equilateral triangles, then there must be a rainbow unit equilateral triangle.
\end{theorem}

Theorem \ref{main} does not work if we have only a constant number of colors. It seems as though one should be able to prove a stronger result, some absolute constant, $c$. However, any such result would have to somehow deal with the relative sizes of the color classes, as the following degenerate example shows.

\begin{example}
Color $\mathbb F_q^2$ by giving each point of the form $(2i,0)$, with $i=1, 2, \dots \lfloor q/2\rfloor,$ a distinct color that isn't blue, and coloring the rest of the points blue. Clearly, none of the pairs of non-blue points are a unit distance apart, so any unit equilateral triangle in this coloring will have at least two blue points.
\end{example}

\section{Proof of Proposition \ref{size2}}

Consider the vertices of any unit equilateral triangle. If they are three different colors, then we are done. If not, call the two vertices that are the same color $a_1$ and $a_2$, and call the third vertex $b$. Note that $b$ is a different color than $a_1$ and $a_2$, as there are no more than two points of any color. Now, there is a unique point, $c_1$, that forms an equilateral triangle with $a_1$ and $b$. Similarly, there is a unique point, $c_2$, that forms an equilateral triangle with $a_2$ and $b$. See Figure \ref{fig1}. It is not possible for both $c_1$ and $c_2$ to be the same color as $b$, so there must exist a rainbow unit equilateral triangle.
\begin{center}\label{fig1}
\includegraphics[scale=1]{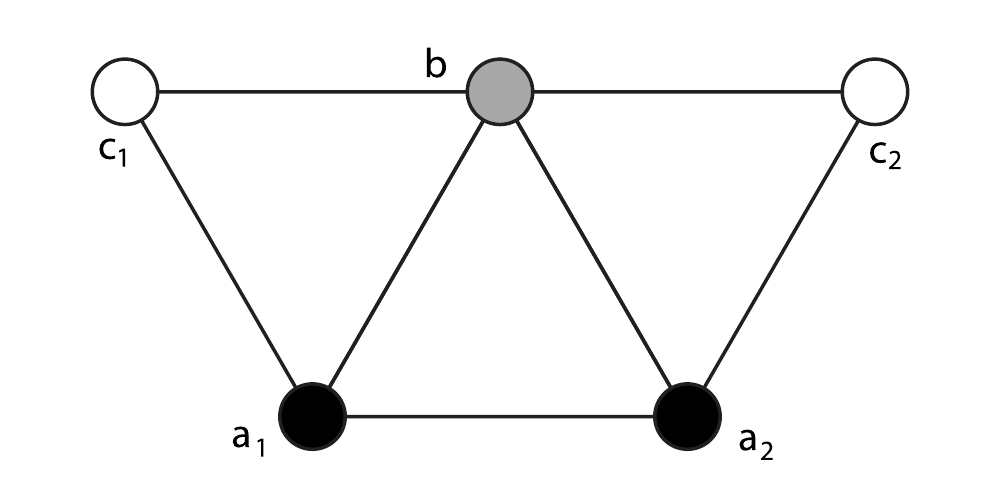}\\
\small{{\bf Figure \ref{fig1}:} The points $a_1$ and $a_2$ are the same color. The point $b$ is a different color, but $c_1$ and $c_2$ cannot both be the same color as $a_1, a_2,$ or $b$.}
\end{center}

\section{Proof of Theorem \ref{main}}

First, there is a minor technical point alluded to in the statement of Theorem \ref{manyColorsTech}, namely that there are vector spaces over finite fields for which no unit equilateral triangles exist. For our purposes, the existence of unit equilateral triangles is equivalent to the presence of an element $s\in \mathbb F_q$ such that $s^2=3.$ This will also guarantee that every pair of points a unit distance apart must determine exactly two unit equilateral triangles, a fact which we often take for granted in the Euclidean setting. We state this precisely as a special case of Lemma 4.1 from \cite{BIP}.

\begin{lemma}\label{lem1}
If there exists an element, $s\in \mathbb F_q$ such that $s^2 = 3$, then for any pair of points, $x,y \in \mathbb F_q^2$, with $d(x,y)=1$, there exist two other points of $\mathbb F_q^2$, $a$ and $b$, such that $x,y,a$ and $x,y,b$ are the vertices of unit equilateral triangles.
\end{lemma}

In what follows, we will always assume that any $q$ under consideration is such that there are equilateral triangles in $\mathbb F_q^2.$

\subsection{Outline}

The outline of the proof is as follows: we begin with two simple counting arguments for constructing a relatively evenly distributed ``sub-coloring" of our original coloring, and then prove a result for colorings where each class roughly the same size. These results are asymptotic in the size of the ambient set.

Define a {\bf fair} $k$-coloring of a set $S$, of size $s$, to be a $k$-coloring where there exists an $n \approx s/k,$ such that each color class has between $.1n$ and $10n$ elements. Notice that this is more relaxed than the related notion of an equinumerable $k$-coloring, where each color class must have size either $\left\lfloor s/k\right\rfloor$ or $\left\lceil s/k\right\rceil$.

Suppose we have a $k$-coloring of a set, $S$. Now, for any $t\geq k$, we will say that a $t$-coloring {\bf refines} the original $k$-coloring if for every color class from the $t$-coloring is contained in a color class from the $k$-coloring. For example, suppose $S:=\{1,2,3,4\}$, and our 2-coloring color classes are $\{1,2\}$ and $\{3,4\}$. Then the 3-coloring $\{1,2\}, \{3\},$ and $\{4\}$ would be a refinement of the 2-coloring, but the 3-coloring $\{1\}, \{2,3\}, \{4\}$ would not be.

\begin{lemma}\label{greedy}
Suppose $S$ is a set with $s$ elements, and we $r$-color it so that no color class has greater than $cs$ elements, for any positive constant $c$. Then there exists a fair $t$-coloring of $S$ for some $t\not\lesssim 1,$ for which the original coloring is a refinement.
\end{lemma}

The next lemma says that if $u\lesssim t$, and we have a fair $t$-coloring of some set, then for some $k\lesssim u$, with $k\not\lesssim 1$, we can use it to generate a fair $k$-coloring of the same set, for which the original coloring is a refinement. Both the statement and the proof are similar to that of Lemma \ref{greedy}, but there are enough subtle differences that we state and prove it separately.

\begin{lemma}\label{manyColors}
Given a fair $t$-coloring of a set $S$ of size $s$, and $u\lesssim s$, with $t,u\not\lesssim 1,$ there exists a fair $k$-coloring of $S$ for some $k\not\lesssim 1,$ with $k \lesssim u$, for which the original coloring is a refinement.
\end{lemma}

The third and final main lemma is a geometric result with some superficial constraints, which are handled by Lemma \ref{greedy} and Lemma \ref{manyColors}.

\begin{lemma}\label{manyColorsTech}
Given a fair $k$-coloring of $\mathbb F_q^2$, with $k \lesssim q^\frac{1}{2}$, and $k\not\lesssim 1,$ if there exist unit equilateral triangles, then there must be a rainbow unit equilateral triangle.
\end{lemma}

We now prove Theorem \ref{main} modulo the proofs of the main lemmas.

\begin{proof}
We begin with a coloring of $\mathbb F_q^2$, where no color class has size greater than $cq^2$, for any constant $c$. So, we apply Lemma \ref{greedy} with $S=\mathbb F_q^2$, and $s=q^2$, and get a fair $t$-coloring of $\mathbb F_q^2,$ for some $t\not\lesssim 1$, for which the original coloring is a refinement. That is, if two points are in different color classes in the new coloring, then they came from different color classes in the original coloring, though the converse may not hold.

Now, we would like to apply Lemma \ref{manyColorsTech} to this new coloring directly, but there is a technical constraint to the methods employed in Lemma \ref{manyColorsTech} that requires us to have rather large color classes. Luckily, Lemma \ref{manyColors} tells us that we can coarsen our coloring again to get a suitable coloring. Specifically, we set $u=q^\frac{1}{2}$, and apply Lemma \ref{manyColors} to our fair $t$-coloring of $\mathbb F_q^2$. 

Applying Lemma \ref{manyColors} in this way guarantees us a fair $k$-coloring of $\mathbb F_q^2$ with $k\lesssim q^\frac{1}{2}$ and $k\not\lesssim 1,$ satisfying the hypotheses of Lemma \ref{manyColorsTech}. As such, we are guaranteed the existence of an equilateral triangle, each of whose vertices comes from a different color class in our newest coloring. Again, this new coloring respects the original coloring in the sense that points from distinct color classes in this new coloring will be in distinct color classes in the original coloring, so this equilateral triangle is rainbow in the original coloring as well.
\end{proof}

\subsection{Proof of Lemma \ref{greedy}}

Essentially, Lemma \ref{greedy} says that if we have a sequence of sets, $S_j$, with each $s_j=|S_j|$, with $s_j>s_{j-1}$, and a sequence of colorings, where we $r_j$-color the set $S_j$, then if the $r_j$ grow with $s_j$, then we can always find a fair $t_j$-coloring of each $S_j$, for some $t_j$ that grows faster than a constant.

\begin{example}
Consider the sequence of the rings of integers of $q$ elements, $\mathbb Z_q$, with $q$ tending to infinity. Suppose that each ring element not equal to 0 or 1 is colored by the largest prime divisor of its corresponding integer, and give 0 and 1 their own color. So the element 6 in each ring will be given the color 3. Letting $q$ tend to infinity, we see that this is a coloring scheme with no color class larger than $c\ln q$, and also with more than $cq$ color classes for any constant $c$. Lemma \ref{greedy} shows that there is a fair $t$-coloring for some $t$ that grows faster than a constant with $q$, for which our original coloring is a refinement. 
\end{example}

\begin{proof}
To prove this, we will check to see if our coloring is already a fair $r$-coloring. If so, set $t=r$, and we are done. If not, we will put two color classes together into a ``super-color," and check again. As the argument is iterative, we will use superscripts on the sets to keep track of which iteration we are on.

Label the color classes $A_j^1$, so that the disjoint union $\dot\cup_j A_j^1 = S,$ and so they are indexed (after a possible relabeling) in non-increasing order,
$$|A_1^1|\geq |A_2^1|\geq \dots \geq |A_f^1|.$$
Set $a=|A_1^1|$, the size of the biggest color class(es). By assumption, $a< cs$ for any positive constant $c$. Now, if all of the color classes are of size between $.1a$ and $10a$, we are done by the definition of a fair coloring.

If not all of the color classes live within this range, then we know that the bound violated is the lower bound, as no set is larger than $a$, by definition. So we take the last two sets, $A_{r-1}^1$ and $A_r^1$ and relabel their union as $A_{r-1}^2$. Now there are two possible outcomes.

If $|A_{r-1}^1| > .9a,$ then because the sets were labeled in non-increasing order, we have that 
\begin{equation}\label{bigA1}
a\geq |A_{j}^1| > .9a, \text{ for all } j=1,2, \dots, r-1.
\end{equation}
Since $|A_{r-1}^1|\leq a$ and $|A_r^1|<.1a$, 
\begin{equation}\label{bigA2}
.9a < |A_{r-1}^2| = |A_{r-1}^1 \cup A_r^1| = |A_{r-1}^1|+|A_r^1|\leq 1.1a.
\end{equation}
Combining \eqref{bigA1} and \eqref{bigA2} give us that our new coloring is fair. Since no set has more than $1.1a$ elements, we have that for any positive constant, $c$, each color class of the new coloring has fewer than $1.1cs$ elements. By the pigeonhole principle, we must have $t$ colors, where $t \not\lesssim s$, and we are done.

Now, if $|A_{r-1}^1|\leq .9a,$ then we know that
$$|A_{r-1}^2| = |A_{r-1}^1 \cup A_r^1| = |A_{r-1}^1|+|A_r^1|\leq a,$$
and $a$ is also the size of the largest color class in the new coloring.

Keep repeating the above procedure. Clearly, the original coloring refines each new coloring, as the new color classes are just unions of old color classes. There are only finitely many color classes, so the process must stop at some point. Notice that the size of the largest color class will remain $a$ throughout each iteration, except possibly the final step, so no color class will ever have a positive proportion of $S$, and our final coloring will have all of the claimed properties.
\end{proof}

\subsection{Proof of Lemma \ref{manyColors}}

\begin{proof}
Notice that if $t\lesssim u$, we set $k=t$, and we are already done. If $t\not\lesssim u$, then we will group color classes together to form $k$ superclasses, for some $k \lesssim u.$

From the definition of fair colorings, there must be a natural number, $m$, such that each color class has size between $.1m$ and $10m$, and $m\approx s/t$. Suppose that the color classes are called $A_j$, where $j=1, \dots, t$, and that they are labeled in non-increasing order, so that
$$|A_1|\geq |A_2| \geq \dots \geq |A_t|.$$

Now, let $\ell = t/u.$ Define $B_1$ to be the union of the sets
$$B_1:= A_1 \cup A_2 \cup \dots \cup A_{\ell_1},$$
where $\ell_1 = c_1\ell$, for some constant $c_2$, so that $|B_1|\leq 10m\ell$, but $|B_1\cup A_{\ell_1+1}|> 10m\ell.$ That is, pack as many of the largest color classes together as we can without the size exceeding $10m\ell.$ Define $B_2$ in a similar manner, as the union of the next $\ell_2$ sets, where $\ell_2=c_2\ell$, for some constant $c_2,$ so that $|B_2|\leq 10m\ell$, but $|B_2\cup A_{\ell_2+1}|> 10m\ell.$

Continue this method of construction for the sets $B_3, B_4,$ and so on. Let $B_{k'}$ be the last such set with $|B_{k'}|\geq .1m\ell.$ Now, if
$$S=\bigcup_{j=1}^{k'} B_j,$$
then we are done, as we set $k=k'$, and the $B_j$ form a fair $k$-coloring meeting all of the requirements above. However, it is possible that there are some sets, $A_{\ell_{k'}+1}, A_{\ell_{k'}+2}, \dots, A_t,$ left over, whose union is strictly less than $.1m\ell.$

If this is the case, we keep the sets $B_1, B_2, \dots, B_{k'-1}$ as before, but go back and define $B_{k'}'$ to be the union of the next $\ell_{k'}'$ sets, where $\ell_{k'}'=c_{k'}'\ell$, for some constant $c_{k'}',$ so that $|B_{k'}'|\leq 7m\ell$, but $|B_{k'}'\cup A_{\ell_{k'}+1}|> 7m\ell.$ Define $B_{k'+1}$ to be the union of the remaining sets. Now, we know that the size of $B_{k'+1}$ is less than $3.1m\ell,$ so we set $k=k'+1$, and the sets $B_1, B_2, \dots B_{k'-1},B_{k'}',$ and $B_{k'+1}$ form a fair $k$-coloring meeting all of the requirements above.
\end{proof}

\subsection{Proof of Lemma \ref{manyColorsTech}}

\begin{proof}
In what follows, we use $C_j$ to indicate various constants, independent of $q$, that we do not compute. Enumerate the color classes, and call each one $A_j$, where $j=1, \dots, k.$ By the definition of a fair $k$-coloring, there exists a natural number, $n\approx q^2/k,$ such that that each color class has between $.1n$ and $10n$ points. To prove Lemma \ref{manyColorsTech}, we will assume that there are no rainbow unit equilateral triangles, and use this assumption to derive a bound on $k$. Therefore, if this bound is violated, then we must be guaranteed the existence of a rainbow unit equilateral triangle.

We also include a special case of Lemma 1.2 from \cite{HIKSU}, which tells us about how many pairs of points are a unit distance apart.

\begin{lemma}\label{lem2}
$$\left|\{(x,y)\in \mathbb F_q^2\times \mathbb F_q^2: d(x,y)=1\}\right|=(1+o(1))q^3.$$
\end{lemma}

Now, Lemma \ref{lem2} tells us that there are $(1+o(1))q^3$ pairs of points that are a unit distance apart, and Lemma \ref{lem1} tells us that each such pair is a part of exactly two unit equilateral triangles. Putting these together with the fact that each triangle has three sides, each of which requires a pair of points, tells us that the number of unit equilateral triangles is approximately $(2/3)$ times the number of unit distances. So we can say that there are $C_1q^3$ unit equilateral triangles in total, regardless of color.

We will suppose that each of these unit equilateral triangles has at least two points from the same color class, which would imply that there are no rainbow unit equilateral triangles in our coloring. Define $u_n$ to be the largest number of pairs of points from any set of $n$ points in $\mathbb F_q^2$, that are a unit distance apart. That is,
$$u_n=\max_{E\subset \mathbb F_q^2; |E|=n}\left|\{(x,y)\in E\times E:d(x,y)=1\}\right|.$$

We now employ the following lemma, which can be found several places, such as \cite{V07} (in the proof of Lemma 2), by Le Anh Vinh.
\begin{lemma}\label{ud}
If $E\subset \mathbb F_q^2$, and $|E|\gtrsim q^\frac{3}{2}$, then the number of unit distances determined by pairs of points from $E$ is $\lesssim q^{-1}|E|^2.$
\end{lemma}

Lemma \ref{ud} tells us that $u_n\lesssim n^2/q$ whenever $n\gtrsim q^\frac{3}{2}.$ Now, define $T$ to be the total number of pairs of points that are both the same color, and separated by a unit distance. Since each of the $k$ color classes can contribute no more than $C_2u_n$ point pairs to $T$, we have that
$$T \lesssim ku_n.$$
Combining this with our bound on $u_n$ yields
$$T\lesssim \frac{kn^2}{q}.$$

But keep in mind that each of our $C_1q^3$ unit equilateral triangles must have a pair of points from the same color class, so we have a lower bound for $T$ as well, namely
$$T\gtrsim q^3.$$
Combining upper and lower bounds on $T$, and estimating $n$ by $q^2/k$ (the possible error here is smaller than the buried constants) tells us
$$q^3 \lesssim \frac{kn^2}{q}\approx \frac{k}{q}\left(\frac{q^2}{k}\right)^2\approx\frac{q^3}{k}$$
$$k \lesssim 1.$$

This tells us that as long as $n\gtrsim q^\frac{3}{2}$ and $k\not\lesssim 1$, we must have a rainbow unit equilateral triangle. As before, we estimate $n$ by $C_3q^2/k$ (the possible error here is, again, buried in the constants) to get the condition that
$$\frac{q^2}{k}\approx n\gtrsim q^\frac{3}{2}$$
$$q^\frac{1}{2}\gtrsim k.$$
\end{proof}

\end{document}